\documentclass[11pt]{amsart}

\usepackage[T1]{fontenc}
\usepackage{lmodern}
\usepackage{amssymb,mathtools}
\usepackage[expansion=false]{microtype}
\usepackage[a4paper,textwidth=13.4cm,textheight=21.3cm,centering]{geometry}
\usepackage{enumitem}
\usepackage{tikz}
\usetikzlibrary{arrows.meta}
\usepackage{needspace}
\usepackage[hidelinks]{hyperref}
\usepackage{bookmark}
\usepackage{placeins}

\setlist[enumerate]{label=\textup{(\roman*)},leftmargin=2.1em,itemsep=3pt}
\numberwithin{equation}{section}

\newtheorem{mainthm}{Theorem}

\newtheorem{theorem}{Theorem}[section]
\newtheorem{proposition}[theorem]{Proposition}
\newtheorem{lemma}[theorem]{Lemma}
\newtheorem{corollary}[theorem]{Corollary}
\theoremstyle{remark}
\newtheorem{remark}[theorem]{Remark}
\newtheorem{example}[theorem]{Example}

\newcommand{\C}{\mathbb C}
\newcommand{\Q}{\mathbb Q}
\newcommand{\Qbar}{\overline{\mathbb Q}}
\newcommand{\Z}{\mathbb Z}
\newcommand{\PP}{\mathbb P}
\newcommand{\W}{\mathcal W}
\newcommand{\MahlerClass}{\Sigma}
\newcommand{\Borel}{\mathcal B}
\newcommand{\DD}{\mathrm D}
\DeclareMathOperator{\conv}{conv}
\DeclareMathOperator{\ord}{ord}
\DeclareMathOperator{\Exp}{Exp}

\DeclareMathOperator{\divisor}{div}

\title[Eremenko's conjecture on Briot--Bouquet equations]{Eremenko's conjecture on higher-order\texorpdfstring{\\}{} Briot--Bouquet equations}
\author{Diego Marques}
\address{Departamento de Matem\'atica, Universidade de Bras\'ilia, Bras\'ilia, DF, Brazil}
\email{diego@mat.unb.br}
\date{}
\hypersetup{
  pdftitle={Eremenko's conjecture on higher-order Briot--Bouquet equations},
  pdfauthor={Diego Marques},
  pdfsubject={Meromorphic differential equations and arithmetic rigidity},
  pdfkeywords={Briot-Bouquet, arithmetic rigidity, Mahler, algebraic independence, indicator diagrams}
}

\begin{document}

\begin{abstract}
The Briot--Bouquet theorem describes the meromorphic solutions of autonomous
first-order algebraic differential equations. We prove Eremenko's conjecture
extending it to all orders: if $f$ is meromorphic in $\mathbb{C}$ and
$P(f,f^{(n)})=0$ for a nonzero polynomial $P$, then $f$ is rational, elliptic,
a rational function of $e^{cz}$ with $c\neq 0$, or satisfies $f^{(n)}=af+b$.
In particular, every meromorphic solution with $P$ irreducible of degree at
least two belongs to the Weierstrass class. Solutions with poles were
classified by Eremenko, Liao and Ng. We settle the remaining entire case with
a rigidity theorem for polynomial compositions: if $h$ is nonconstant entire
and $D^nR(h)=S(h)$, where $D=d/dz$ and $\deg R=\deg S\ge 2$, then
$h(z)=Ae^{bz}+Be^{-bz}+C$. No growth assumption is needed. The proof turns
the cancellation of leading terms into a contraction of indicator diagrams
and obtains incompatible zero counts near an exposed edge. As an arithmetic
application, an algebraic Laurent expansion at one algebraic point and an
algebraic value at another force a solution to lie in
$\overline{\mathbb{Q}}(z)$. Hence any two distinct derivatives of a function
in Mahler's class (transcendental meromorphic functions with algebraic
Laurent coefficients at every algebraic point) are algebraically independent
over $\mathbb{C}$.
\end{abstract}

\maketitle

\section{Introduction}\label{sec:intro}

The classical Briot--Bouquet theorem classifies the meromorphic
solutions on $\C$ of an autonomous first-order algebraic differential
equation as elliptic functions and their degenerations. In higher
order, this conclusion fails already for $f'''=f$. Eremenko's
conjecture asserts that affine equations account for all exceptions.
We prove this conjecture by resolving its remaining entire case.
The main new ingredient is a rigidity theorem for nonlinear
polynomial compositions. As an arithmetic consequence, any two
distinct derivative levels of a function in Mahler's class are
algebraically independent over $\C$.

\Needspace{8\baselineskip}
\subsection{The classification and its analytic core}
Throughout, meromorphic functions are single-valued on $\C$ and are
not identically infinite. A meromorphic function is transcendental if
it is not rational. The \emph{Weierstrass class} $\W$ consists of
rational functions, elliptic functions, and functions $R(e^{az})$ with
$R\in\C(T)$ and $a\neq0$. We consider
\begin{equation}\label{eq:BB}
 P\bigl(f,f^{(n)}\bigr)=0,
 \qquad 0\neq P\in\C[X,Y],\quad n\geq1.
\end{equation}
Polynomial degrees are total degrees unless stated otherwise.

\begin{mainthm}[Eremenko's conjecture]\label{thm:classification}
Let $n\geq1$, and let $f$ be meromorphic on $\C$. If $f$ and $f^{(n)}$
are algebraically dependent over $\C$, then either $f\in\W$, or
\begin{equation}\label{eq:affine}
 f^{(n)}=af+b\qquad\text{for some }a,b\in\C.
\end{equation}
In particular, every meromorphic solution of an irreducible equation
\eqref{eq:BB} with $\deg P\geq2$ belongs to $\W$.
\end{mainthm}

For $n\leq2$, the affine solutions already belong to $\W$;
exceptions outside $\W$ can therefore occur only for $n\geq3$.
Every entire solution is an exponential polynomial: this follows
from the description of entire members of $\W$ and the explicit
solutions of \eqref{eq:affine}. The analytic core is the following
statement, where $\DD=d/dz$.

\begin{mainthm}[Polynomial-composition rigidity]\label{thm:composition}
Let $n\geq1$ and let $R,S\in\C[T]$ have the same degree $d\geq2$.
Suppose that $h$ is nonconstant and entire and
\begin{equation}\label{eq:RS}
 \DD^nR(h)=S(h).
\end{equation}
Then
\begin{equation}\label{eq:h-classified}
 h(z)=Ae^{bz}+Be^{-bz}+C
\end{equation}
for some $A,B,C,b\in\C$ with $b\neq0$ and $(A,B)\neq(0,0)$.
In particular, $R(h)$ is a Laurent polynomial in an exponential.
\end{mainthm}

Finite positive exponential type follows from the equation itself;
see Lemma~\ref{lem:automatic-growth}.

\noindent\emph{The degree restriction is essential.}
If $\omega^3=1$ and $\omega\neq1$, then
\begin{equation}\label{eq:degree-one-example}
 h(z)=e^z+e^{\omega z}+e^{\omega^2z}
\end{equation}
satisfies $h'''=h$ and has finite positive exponential type, but is not
of the form \eqref{eq:h-classified}. Thus Theorem~\ref{thm:composition}
fails for $R(T)=S(T)=T$. The same example explains the affine
alternative in Theorem~\ref{thm:classification}.

\subsection{The remaining case}
The first-order theorem goes back to Briot and Bouquet
\cite{BriotBouquet}. Picard \cite{Picard} proved the second-order
case in 1880, later conjectured again by Hille \cite{Hille} and
reproved by Bank and Kaufman \cite{BankKaufman}. Eremenko developed
the higher-order theory in \cite{Eremenko82,EremenkoSurvey};
Eremenko, Liao and Ng \cite[Theorem~1]{ELN} subsequently classified
all solutions having a pole.

For entire solutions, Yanchenko \cite{Yanchenko} treated the case
of distinct leading linear factors. His 2024 theorem
\cite[Russian original, p.~1132]{Yanchenko24} allows lower-degree terms depending polynomially
on $z$ and on additional derivative levels, but requires the leading
part $\prod_{j=1}^{d}(f^{(\ell)}-\alpha_jf)$ to have pairwise
distinct factors; the remainder has degree at most $d-1$ in the
dependent variables. As clarified in Eremenko's December~2024
addendum \cite{EremenkoNote}, the results of \cite{Eremenko82}
already imply that the highest homogeneous part of an irreducible
Briot--Bouquet equation admitting a transcendental entire solution
has at most two distinct linear factors, and they settle the
two-factor case. The remaining configuration is
\begin{equation}\label{eq:degenerate-case}
 \bigl(f^{(n)}-af\bigr)^d+Q\bigl(f,f^{(n)}\bigr)=0,
 \qquad d\geq2,\quad\deg Q<d.
\end{equation}
For example, $f(z)=e^{2z}+e^z$ satisfies
\[
 (f'''-8f)^2-7f'''+7f=0.
\]
The defining polynomial is irreducible: with $T=Y-8X$, it becomes
$T^2-7T-49X$. Its highest homogeneous part is $(Y-8X)^2$.
Thus the repeated-factor case occurs already in a quadratic equation.

A 2026 correction to Merolla's article
\cite{MerollaCorrection} withdraws an earlier claimed proof
of the conjecture because it relied on an incorrect auxiliary theorem.

For transcendental entire solutions of \eqref{eq:BB}, the
necessary growth estimate was already available: Theorem~D
and Lemma~1 of \cite{ELN} imply order one and finite positive
exponential type (see Proposition~\ref{prop:known-input}).
The remaining difficulty was to determine the form of the entire
parametrizing function. Indeed, normalization and Picard's theorem reduce the case outside
$\W$ to
\[
 f=R(h),\qquad f^{(n)}=S(h),\qquad \deg R=\deg S=d,
\]
with $h$ entire of finite positive exponential type. For $d=1$ the
relation is affine; for $d\geq2$, Theorem~\ref{thm:composition}
completes the classification. The one-point-at-infinity reduction
has precisely the repeated leading factor in
\eqref{eq:degenerate-case}; see Remark~\ref{rem:infinity}.

To indicate the mechanism, write $a=s_d/r_d$ for the ratio of the
leading coefficients of $S$ and $R$, and put $F=R(h)$ and
$Q=S-aR$. Then
\begin{equation}\label{eq:composition-intro}
 (\DD^n-a)F=Q(h),\qquad\deg Q<d.
\end{equation}
The degree drop contracts the augmented indicator diagram of the
right-hand side by a factor smaller than one. Borel inversion
separates $F$ into modes $e^{\lambda z}$, with $\lambda^n=a$, and a
remainder supported in the contracted diagram. A non-antipodal
exposed edge produces discs $D_k$ on which
\[
 N_{D_k}(F)=1,\qquad N_{D_k}(F)=dN_{D_k}(h).
\]
The first count follows by comparison with two exponential modes;
the second uses polynomial properness and Rouch\'e's theorem.
This is where $d\geq2$ is decisive. The surviving spectrum consists
of one frequency or an antipodal pair. Sectorial inversion then
leaves a remainder of exponential type zero, which the differential
equation forces to be constant. No regular-growth assumption is
made in this argument.

\subsection{Arithmetic rigidity and Mahler's question}
Let $\Qbar$ be the algebraic closure of $\Q$ in $\C$. A pole is
regarded as the projective value $\infty\in\PP^1(\Qbar)$.
The classification has the following consequence.

\begin{mainthm}[Two-point arithmetic rigidity]\label{thm:germ}
Let $f$ be meromorphic on $\C$, and let $\alpha,\beta\in\Qbar$ be
distinct. Suppose that
\[
 f(\alpha+t)\in\Qbar((t)),\qquad f(\beta)\in\PP^1(\Qbar).
\]
If $f$ and $f^{(n)}$ are algebraically dependent over $\C$ for some
$n\geq1$, then $f\in\Qbar(z)$.
\end{mainthm}

The first hypothesis concerns the full convergent Laurent expansion,
not a finite jet or an algebraic function germ. It descends a finite
differential relation to $\Qbar$. The second hypothesis prescribes
only one value, not its derivatives. After classification, classical
transcendence theorems exclude the exponential and elliptic
possibilities. Consequently, if $f$ is transcendental, $\alpha\in\Qbar$,
$f(\alpha+t)\in\Qbar((t))$, and $f$ and $f^{(n)}$ are
algebraically dependent over $\C$ for some $n\geq1$, then
\begin{equation}\label{eq:exceptional-set}
 \{\beta\in\Qbar:f(\beta)\in\PP^1(\Qbar)\}=\{\alpha\}.
\end{equation}

In \cite[Section~5, p.~104]{Mahler}, Mahler considers the class
$\MahlerClass$ of transcendental meromorphic functions with algebraic
Laurent coefficients at every algebraic point, and asks:
\begin{quote}
\itshape
Can a function $f(z)$ in $\Sigma$ satisfy an algebraic differential
equation?
\end{quote}
Thus $f\in\MahlerClass$ means that, for every $\alpha\in\Qbar$,
\begin{equation}\label{eq:mahler-class}
 f(z)=\sum_{j\geq-N_\alpha}a_{\alpha,j}(z-\alpha)^j,
 \qquad a_{\alpha,j}\in\Qbar.
\end{equation}
For entire $f$, this is equivalent to
\begin{equation}\label{eq:entire-class}
 f^{(j)}(\Qbar)\subseteq\Qbar\qquad(j\geq0).
\end{equation}

Faber's construction shows that $\MahlerClass$ contains entire
functions; see \cite[Section~5]{Mahler}. More generally,
Huang, Marques and Mereb \cite{HMM} prescribe dense target sets for
all derivatives simultaneously on any countable set of source
points. Such interpolation is therefore flexible without a
differential constraint. Neither a common number field nor uniform
arithmetic bounds are imposed in the definition of $\MahlerClass$.

\begin{corollary}\label{cor:allpairs}
Let $f\in\MahlerClass$ and let $0\leq r<s$. Then $f^{(r)}$ and
$f^{(s)}$ are algebraically independent over $\C$.
\end{corollary}

Only two algebraic source points are used in this application,
although membership in $\MahlerClass$ prescribes all of them.
The conclusion is independence over the \emph{constant field} $\C$.
It does not exclude relations
\begin{equation}\label{eq:full-mahler}
 P\bigl(z,f,f',\ldots,f^{(m)}\bigr)=0,
 \qquad0\neq P\in\C[Z,Y_0,\ldots,Y_m],
\end{equation}
with explicit dependence on $z$, or autonomous relations involving
three or more derivative levels. The full question of differential
transcendence remains outside the present result.

The arithmetic step uses the classical transcendence theorems
of Hermite--Lindemann, Lindemann--Weierstrass and Schneider
\cite{Schneider,Lang,Waldschmidt}. The Schneider--Lang theorem
does not directly apply: Mahler's condition provides neither
simultaneous values in a fixed number field nor a general growth
bound. Here the classification reduces the argument to explicit
transcendence obstructions at two points.

\subsection{Related rigidity results and organization}
Related rationality alternatives arise from simultaneous
differential and difference constraints \cite{ADH}, or from two
sufficiently independent difference operators \cite{ADHW}.
Here the constraint is an autonomous, possibly nonlinear
differential relation; global classification and arithmetic descent
replace the additional difference structure.

Section~\ref{sec:borel} records the indicator-diagram and
Borel-transform tools. Section~\ref{sec:composition} proves
Theorem~\ref{thm:composition}, and
Section~\ref{sec:classification} deduces
Theorem~\ref{thm:classification}. Sections~\ref{sec:arithmetic}
and \ref{sec:main-proof} establish the arithmetic results.
Sections~\ref{sec:scope} and \ref{sec:further} discuss sharpness,
obstructions to extensions, and further questions. Zeros are counted
with multiplicity throughout.

\section{Indicator diagrams and Borel inversion}\label{sec:borel}

For a nonzero entire function $u$, put
\[
 M(r,u)=\max_{|z|=r}|u(z)|,\qquad
 \tau(u)=\limsup_{r\to\infty}\frac{\log M(r,u)}r,
\]
and set $\tau(0)=0$.
We write $\log^+x=\max\{\log x,0\}$ for $x>0$, with $\log^+0=0$.
We say that $u$ has exponential type if $\tau(u)<\infty$, and type zero if $\tau(u)=0$.
For an entire function of exponential type, the indicator
\[
 h_u(\theta)=\limsup_{r\to\infty}\frac{\log|u(re^{i\theta})|}{r}
\]
is the support function of its conjugate indicator diagram $I(u)$. We use the convention
\begin{equation}\label{eq:support}
 H_L(z)=\max_{\zeta\in L}\operatorname{Re}(\zeta z),
 \qquad H_{I(u)}(re^{i\theta})=r h_u(\theta).
\end{equation}
Our convention is the real bilinear pairing $\operatorname{Re}(\zeta z)$,
not the Hermitian pairing $\operatorname{Re}(\zeta\overline z)$.
Accordingly, $I(e^{\lambda z})=\{\lambda\}$, and the maximizing
unit direction for a nonzero frequency $\lambda$ is
$\overline\lambda/|\lambda|$. Set $I(0)=\varnothing$ and define the augmented diagram
\[
 K(u)=\conv\bigl(I(u)\cup\{0\}\bigr),\qquad K(0)=\{0\}.
\]
We use the conventions $\log 0=-\infty$ and $h_0=-\infty$;
support-function identities involving $I(u)$ are understood for
$u\neq0$.

Then
\begin{equation}\label{eq:type-diagram}
 H_{K(u)}(e^{i\theta})=\max\{h_u(\theta),0\},
 \qquad \tau(u)=\max_{\zeta\in K(u)}|\zeta|.
\end{equation}
In particular, $u$ has type zero if and only if $K(u)=\{0\}$. For $c\neq0$, change of variables in the indicator gives
\begin{equation}\label{eq:scaling-diagram}
 I\bigl(u(c\,\cdot)\bigr)=cI(u),\qquad
 K\bigl(u(c\,\cdot)\bigr)=cK(u).
\end{equation}

We recall the Borel--P\'olya representation and the associated support
estimates; see \cite[Section~2, Proposition~2.1]{BeiseMuller}
and the classical treatments \cite{Boas,Levin}.
For a nonempty compact convex set $L\subset\C$, write $u\in\Exp(L)$ if $u$ is entire of exponential type and $I(u)\subseteq L$.
The zero function belongs to every $\Exp(L)$.
Equivalently, for every $\varepsilon>0$ there is $C_\varepsilon>0$ such that
\begin{equation}\label{eq:uniform-support}
 |u(z)|\leq C_\varepsilon
 \exp\bigl(H_L(z)+\varepsilon|z|\bigr)
 \qquad(z\in\C).
\end{equation}
The Borel transform
\begin{equation}\label{eq:borel}
 \Borel u(\zeta)=\sum_{j\geq0}\frac{u^{(j)}(0)}{\zeta^{j+1}}
\end{equation}
converges near infinity, vanishes there, and extends holomorphically to $\C\setminus I(u)$.
If $\Gamma$ is a positively oriented simple contour enclosing $I(u)$ and lying in its complement, then
\begin{equation}\label{eq:polya}
 u(z)=\frac1{2\pi i}\int_\Gamma\Borel u(\zeta)e^{\zeta z}\,d\zeta.
\end{equation}
The contour may be chosen arbitrarily close to any compact convex set containing $I(u)$.
The indicator of a nonzero function of exponential type is continuous, since it is a support function of a compact convex set.
Estimate~\eqref{eq:uniform-support}, rather than a pointwise directional estimate, will control the moving circles used in Section~\ref{sec:composition}.

\begin{lemma}[Polynomial composition]\label{lem:poly-diagram}
If $u$ is entire of exponential type and $R\in\C[T]$ has degree $d\geq1$, then
\[
 K(R(u))=dK(u),\qquad \tau(R(u))=d\tau(u).
\]
If $R$ is constant, then $K(R(u))=\{0\}$.
\end{lemma}
\begin{proof}
The leading term of $R$ gives
\[
 \log^+|R(w)|=d\log^+|w|+O_R(1)
 \qquad(w\in\C),
\]
where the bound is uniform: it follows from leading-term comparison for large $|w|$ and boundedness on a fixed disc for the remaining values.
Substitute $w=u(re^{i\theta})$, divide by $r$, and pass to the limsup.
The resulting identity is
\[
 \max\{h_{R(u)}(\theta),0\}
 =d\max\{h_u(\theta),0\}.
\]
These are the support functions of $K(R(u))$ and $dK(u)$, respectively.
Equality of the diagrams and \eqref{eq:type-diagram} give the assertions.
The constant case is immediate.
\end{proof}

\begin{lemma}[Inversion away from the spectrum]\label{lem:inverse}
Let $L\subset\C$ be nonempty, compact, and convex.
Suppose that $G\in\Exp(L)$ and that $p\in\C[\zeta]$ has no zero on $L$.
There exists $v\in\Exp(L)$ such that $p(\DD)v=G$.
\end{lemma}
\begin{proof}
Since the zero set of $p$ is disjoint from the compact set $L$, a sufficiently small neighbourhood of $L$ contains no zero of $p$.
Choose a positively oriented contour $\Gamma$ around $L$ in that neighbourhood, with its interior also avoiding the zeros of $p$, and put
\[
 v(z)=\frac1{2\pi i}\int_\Gamma
       \frac{\Borel G(\zeta)}{p(\zeta)}e^{\zeta z}\,d\zeta.
\]
The integrand is holomorphic off $L$ in this neighbourhood. Any two
such contours have the same winding number about $L$ and can be
compared by Cauchy's theorem, so the integral is independent of the
contour. Differentiation under the integral, followed by \eqref{eq:polya}, gives $p(\DD)v=G$.

For every sufficiently small $\varepsilon>0$, choose $\Gamma$ in the $\varepsilon$-neighbourhood of $L$. Each $\zeta\in\Gamma$ then satisfies
\[
 \operatorname{Re}(\zeta z)\leq H_L(z)+\varepsilon|z|.
\]
The contour integral yields \eqref{eq:uniform-support}, with a constant depending on the contour and hence on $\varepsilon$.
Contour independence ensures that all these estimates concern the same function $v$.
Thus $v\in\Exp(L)$.
\end{proof}

\begin{lemma}[Separation over type-zero coefficients]\label{lem:zero-independence}
Let $\mu_1,\ldots,\mu_s$ be distinct complex numbers, and let $u_1,\ldots,u_s$ be entire functions of exponential type zero. If
\[
 \sum_{j=1}^s u_j(z)e^{\mu_jz}=0,
\]
then $u_1=\cdots=u_s=0$.
\end{lemma}
\begin{proof}
The transform $\Borel u_j$ is holomorphic on $\C\setminus\{0\}$ and vanishes at infinity.
The product rule, applied in a neighbourhood of infinity where the series converge absolutely, gives
\begin{equation}\label{eq:borel-shift}
 \Borel\bigl(u_je^{\mu_j\,\cdot}\bigr)(\zeta)
 =\Borel u_j(\zeta-\mu_j).
\end{equation}
Indeed, for each $\ell\geq0$,
\[
 \sum_{m\geq\ell}\binom m\ell
 \frac{\mu^{m-\ell}}{\zeta^{m+1}}
 =\frac1{(\zeta-\mu)^{\ell+1}}.
\]
Taking Borel transforms of the given identity and using analytic continuation yields
\[
 \sum_{j=1}^s\Borel u_j(\zeta-\mu_j)=0
 \qquad\bigl(\zeta\notin\{\mu_1,\ldots,\mu_s\}\bigr).
\]
Near $\mu_j$, every summand other than the $j$th is holomorphic.
Hence the possible isolated singularity of $\Borel u_j(\zeta-\mu_j)$ at $\mu_j$ is removable.
Thus $\Borel u_j$ is entire and vanishes at infinity.
Liouville's theorem makes it zero, and \eqref{eq:polya} gives $u_j=0$.
\end{proof}

Sums, products, and derivatives of type-zero entire functions again have type zero.
For sums and products, this follows from bounds of the form $C_\varepsilon e^{\varepsilon|z|}$; for derivatives, use Cauchy's estimate on unit circles.
We will also need the following equality, including its reverse inequality.

\begin{lemma}[Differentiation preserves exponential type]\label{lem:derivative-type}
If $u$ is entire of exponential type and $m\geq1$, then
\[
 \tau\bigl(u^{(m)}\bigr)=\tau(u).
\]
\end{lemma}
\begin{proof}
Cauchy's estimate on unit circles gives
\[
 M(r,u^{(m)})\leq m!M(r+1,u),
\]
so $\tau(u^{(m)})\leq\tau(u)$. Conversely, Taylor's formula with integral remainder on the segment from $0$ to $z$ gives
\[
 u(z)=\sum_{j=0}^{m-1}\frac{u^{(j)}(0)}{j!}z^j
 +\frac{z^m}{(m-1)!}\int_0^1(1-s)^{m-1}u^{(m)}(sz)\,ds.
\]
Consequently,
\[
 M(r,u)\leq C(1+r^{m-1})+\frac{r^m}{m!}M(r,u^{(m)}).
\]
This proves the reverse inequality. If $u^{(m)}=0$, then $u$ is a polynomial and both types are zero.
\end{proof}

\begin{lemma}[Automatic growth in the composition equation]
\label{lem:automatic-growth}
Let $n\geq1$, let $R,S\in\C[T]$ have the same degree $d\geq1$,
and let $h$ be a nonconstant entire function satisfying
$\DD^nR(h)=S(h)$. Then $h$ and $R(h)$ have finite positive
exponential type.
\end{lemma}
\begin{proof}
Let $r_d,s_d$ be the leading coefficients of $R,S$, and put
$F=R(h)$, $a=s_d/r_d\neq0$, and $Q=S-aR$. Then $F$ is nonconstant and $\deg Q<d$.
Polynomial comparison gives constants $C>0$ and $0\leq\gamma<1$
such that, on $\C$,
\[
 |F^{(n)}-aF|=|Q(h)|\leq C(1+|F|)^\gamma.
\]
For nonconstant $Q$ one can take $\gamma=\deg Q/d$; for constant
or zero $Q$ take $\gamma=0$.
In particular, $|F^{(n)}|\leq C_1(1+|F|)$. Along each ray,
\[
 U_\theta(t)=1+\sum_{j=0}^{n-1}|F^{(j)}(te^{i\theta})|
\]
is locally absolutely continuous and satisfies
$U_\theta'\leq C_2U_\theta$ almost everywhere, with $C_2$ and
$U_\theta(0)$ independent of $\theta$. Gronwall's inequality gives
$M(r,F)\leq C_3e^{C_2r}$.

Choose $z_r$ with $|z_r|=r$ and $|F(z_r)|=M(r,F)$.
Since $F$ is nonconstant, $M(r,F)\to\infty$, so the first estimate
implies $|F^{(n)}(z_r)|\geq |a|M(r,F)/2$ for all sufficiently large $r$.
Fix $L>0$ with $\eta:=|a|L^n/(2n!)>1$. Cauchy's estimate on
$D(z_r,L)$ then yields
\[
 M(r+L,F)\geq\frac{L^n}{n!}|F^{(n)}(z_r)|
             \geq\eta M(r,F).
\]
Iteration gives $\tau(F)\geq L^{-1}\log\eta>0$.
Finally, polynomial comparison gives
$|h|^d\leq C_4(1+|F|)$, so $h$ has finite exponential type.
Lemma~\ref{lem:poly-diagram} now yields $\tau(F)=d\tau(h)$,
and hence $0<\tau(h)<\infty$.
\end{proof}

\section{Polynomial-composition rigidity}\label{sec:composition}

We prove Theorem~\ref{thm:composition}. The equality of degrees yields
a strict contraction after cancellation of the leading terms. The
hypothesis $d\geq2$ then enters through the incompatibility of two
local zero counts.

Write
\[
 R(T)=r_dT^d+\cdots+r_0,\qquad
 S(T)=s_dT^d+\cdots+s_0,
\]
and put
\begin{equation}\label{eq:F-G}
 a=\frac{s_d}{r_d}\neq0,\qquad
 F=R(h),\qquad Q=S-aR,\qquad G=Q(h).
\end{equation}
Then
\[
 (\DD^n-a)F=G,\qquad \deg Q<d.
\]
By Lemma~\ref{lem:automatic-growth}, $h$ and $F$ have finite
positive exponential type.
Set $m_Q=\deg Q$ when $Q$ is nonconstant, and $m_Q=0$ when $Q$ is constant or zero. Define
\[
 \vartheta=\frac{m_Q}{d}<1,\qquad \rho_*=|a|^{1/n},\qquad K=K(F).
\]
Lemma~\ref{lem:poly-diagram} gives the exact contraction identities
\begin{equation}\label{eq:contraction}
 K(G)=\vartheta K,\qquad \tau(G)=\vartheta\tau(F).
\end{equation}
Notice that $\tau(F)=d\tau(h)>0$.

\subsection{The contracted remainder and the frequency polygon}

\begin{lemma}\label{lem:polygon}
There are a nonempty set $\Lambda\subseteq\{\lambda:\lambda^n=a\}$ and nonzero constants $c_\lambda$ such that
\begin{equation}\label{eq:spectral}
 F(z)=\sum_{\lambda\in\Lambda}c_\lambda e^{\lambda z}+v(z),
 \qquad v\in\Exp(\vartheta K).
\end{equation}
Moreover,
\begin{equation}\label{eq:polygon}
 K=\conv\bigl(\{0\}\cup\Lambda\bigr),\qquad \tau(F)=\rho_*.
\end{equation}
\end{lemma}
\begin{proof}
Let $p(\zeta)=\zeta^n-a$. Its roots are simple and lie on $|\zeta|=\rho_*$. Variation of constants gives
\begin{equation}\label{eq:variation}
 F(z)=\sum_{\lambda^n=a}b_\lambda e^{\lambda z}
 +\sum_{\lambda^n=a}\frac1{p'(\lambda)}
   \int_0^z e^{\lambda(z-t)}G(t)\,dt,
\end{equation}
where the constants account for the initial data. Integrate on $t=sz$, $0\leq s\leq1$. For every $\varepsilon>0$, the exponential-type bound on $G$ gives
\[
 \left|\int_0^z e^{\lambda(z-t)}G(t)\,dt\right|
 \leq C_\varepsilon |z|
 e^{\max\{\rho_*,\tau(G)+\varepsilon\}|z|}.
\]
It follows that
\begin{equation}\label{eq:initial-type-bound}
 \tau(F)\leq\max\{\rho_*,\tau(G)\}.
\end{equation}
Together with \eqref{eq:contraction} and $\tau(F)>0$, this implies $\tau(F)\leq \rho_*$: otherwise both $\rho_*$ and $\vartheta\tau(F)$ would be smaller than $\tau(F)$.

Thus $K$ lies in the closed disc of radius $\rho_*$, and $\vartheta K$ lies in the smaller disc of radius $\vartheta\rho_*$.
In particular, $p$ has no zero on $\vartheta K$.
Lemma~\ref{lem:inverse} provides a particular solution $v\in\Exp(\vartheta K)$ of $p(\DD)v=G$.
The difference $F-v$ belongs to $\ker(\DD^n-a)$, which gives \eqref{eq:spectral} after discarding the zero coefficients.
If none remained, then $K\subseteq \vartheta K$, forcing $K=\{0\}$ by boundedness and $\vartheta<1$.
This contradicts $\tau(F)>0$, so $\Lambda$ is nonempty.

Fix $\lambda\in\Lambda$ and take the unit direction
\[
 \nu=\frac{\overline\lambda}{\rho_*}.
\]
This conjugation is determined by the pairing in \eqref{eq:support}. We have
\[
 \operatorname{Re}(\lambda\nu)=\rho_*,
 \qquad \operatorname{Re}(\mu\nu)<\rho_*
 \quad(\mu\in\Lambda\setminus\{\lambda\}).
\]
Since $v$ has type at most $\vartheta\rho_*<\rho_*$, the $\lambda$-mode strictly dominates every other term on $z=t\nu$ as $t\to+\infty$.
Hence $h_F(\arg\nu)=\rho_*$.
The diagram $K$ is contained in $|\zeta|\leq \rho_*$, and the unique point of this disc satisfying $\operatorname{Re}(\zeta\nu)=\rho_*$ is $\lambda$.
Thus $\lambda\in K$.

Put $K_0=\conv(\{0\}\cup\Lambda)$. We have proved $K_0\subseteq K$. Conversely, the indicator bound for a finite sum in \eqref{eq:spectral} yields
\[
 H_K(z)\leq\max\{H_{K_0}(z),\vartheta H_K(z)\}
 \qquad(z\in\C).
\]
Both support functions are nonnegative.
If $H_K(z)>H_{K_0}(z)$ at some $z$, the last inequality would give $H_K(z)\leq \vartheta H_K(z)$ with $H_K(z)>0$, a contradiction.
Therefore $K\subseteq K_0$, proving the diagram identity.
Since $\Lambda$ contains a point of modulus $\rho_*$, the type identity follows as well.
\end{proof}

\subsection{An exposed-edge obstruction}

The next argument is where the nonlinear outer polynomial is used.
For a disc $D$ and a holomorphic function $g$ without zeros on $\partial D$, write $N_D(g)$ for the number of zeros of $g$ in $D$, counted with multiplicity.

\begin{lemma}\label{lem:edge}
The set $\Lambda$ in Lemma~\ref{lem:polygon} consists of one point or two antipodal points.
\end{lemma}
\begin{proof}
Suppose otherwise.
The circle $|\lambda|=\rho_*$ has two consecutive points of $\Lambda$ bounding an empty arc of angular length less than $\pi$.
For at least three points, this follows because the circular gaps sum to $2\pi$; for two non-antipodal points, choose the shorter arc.
We may write the endpoints as
\[
 \lambda_1=\rho_*e^{i(\phi-\alpha)},\qquad
 \lambda_2=\rho_*e^{i(\phi+\alpha)},\qquad
 0<\alpha<\frac\pi2.
\]
With $\nu=e^{-i\phi}$ and $\sigma_*=\rho_*\cos\alpha$, we obtain
\begin{equation}\label{eq:positive-edge}
 \begin{split}
 \operatorname{Re}(\lambda_1\nu)
 &=\operatorname{Re}(\lambda_2\nu)=\sigma_*>0,\\
 \operatorname{Re}(\lambda\nu)&<\sigma_*
 \qquad(\lambda\in\Lambda\setminus\{\lambda_1,\lambda_2\}).
 \end{split}
\end{equation}
Thus the chord is an exposed edge of $K=\conv(\{0\}\cup\Lambda)$; adjoining the origin does not change the edge because $\sigma_*>0$. Figure~\ref{fig:edge} illustrates the two planes used in the argument.

Set $c_j=c_{\lambda_j}$ and
\[
 E(z)=c_1e^{\lambda_1z}+c_2e^{\lambda_2z}.
\]
Its zeros are
\begin{equation}\label{eq:edge-zeros}
 z_k=\frac{\ell+2\pi i k}{\lambda_1-\lambda_2},\qquad
 k\in\Z,\qquad e^\ell=-\frac{c_2}{c_1}.
\end{equation}
They are simple, since
\[
 E'(z_k)=c_1e^{\lambda_1z_k}(\lambda_1-\lambda_2)\neq0.
\]
Indeed,
\[
\frac{2\pi i}{\lambda_1-\lambda_2}
=-\frac{\pi}{\rho_*\sin\alpha}\,\nu.
\]
Choosing the indices tending to $-\infty$ in
\eqref{eq:edge-zeros} and relabelling, we obtain
\[
 z_k=t_k\nu+z_*,\qquad t_k\longrightarrow+\infty,
\]
with $z_*$ fixed.

Choose $\delta>0$ smaller than half the distance between consecutive zeros in \eqref{eq:edge-zeros}, and let $D_k=D(z_k,\delta)$. For $|w|=\delta$,
\[
 E(z_k+w)=c_1e^{\lambda_1z_k}
 \bigl(e^{\lambda_1w}-e^{\lambda_2w}\bigr).
\]
The factor in parentheses has no zero on $|w|=\delta$, so its minimum modulus there is positive. Since
$|e^{\lambda_1z_k}|=e^{\sigma_*t_k+\operatorname{Re}(\lambda_1z_*)}$, there is $c>0$, independent of $k$, such that
\begin{equation}\label{eq:edge-large}
 \min_{\partial D_k}|E|\geq ce^{\sigma_*t_k}.
\end{equation}

By \eqref{eq:positive-edge}, the remaining exponential modes are $O(e^{(\sigma_*-\eta)t_k})$ uniformly on $\partial D_k$, for some fixed $\eta>0$ when such modes are present.
For the remainder, boundedness of $K$ and $z-z_k$ gives
\[
 H_K(z)\leq t_kH_K(\nu)+O(1)=\sigma_*t_k+O(1),
 \qquad |z|\leq t_k+O(1)
 \quad(z\in\partial D_k).
\]
Estimate~\eqref{eq:uniform-support} therefore yields
\begin{equation}\label{eq:edge-remainder}
 \max_{\partial D_k}|v|
 \leq C_\varepsilon\exp\bigl(\vartheta\sigma_*t_k+\varepsilon t_k+O(1)\bigr).
\end{equation}
Take $\varepsilon=(1-\vartheta)\sigma_*/2$. Equations~\eqref{eq:edge-large}
and \eqref{eq:edge-remainder}, together with the estimate for the other
modes, give constants $C_1,\gamma>0$, independent of $k$, such that
\[
 \max_{\partial D_k}\left|\frac{F-E}{E}\right|
 \leq C_1 e^{-\gamma t_k}.
\]
For all large $k$, Rouch\'e's theorem gives
\begin{equation}\label{eq:first-zero-count}
 N_{D_k}(F)=N_{D_k}(E)=1.
\end{equation}
The same estimates give the separate boundary conclusion
\begin{equation}\label{eq:boundary-properness}
 \min_{\partial D_k}|F|\geq\frac c2e^{\sigma_*t_k}
 \longrightarrow\infty.
\end{equation}

Now use $F=R(h)$. For every $M>0$, the quantity
\[
 C_M=\max_{|w|\leq M}|R(w)|
\]
is finite. Once $\min_{\partial D_k}|F|>C_M$, no value of $h$ on $\partial D_k$ can have modulus at most $M$. Thus
\[
 \min_{\partial D_k}|h|\longrightarrow\infty.
\]
Choose $M$ large enough that
\[
 |R(w)-r_dw^d|<|r_dw^d|\qquad(|w|>M).
\]
For all sufficiently large $k$, this inequality holds with $w=h(z)$ throughout $\partial D_k$. A second application of Rouch\'e's theorem gives
\begin{equation}\label{eq:second-zero-count}
 N_{D_k}(F)=N_{D_k}(r_dh^d)=dN_{D_k}(h).
\end{equation}
This compares total zero counts; it does not identify the zeros of $R(h)$ with those of $h$.
In particular, nonzero roots and the constant term of $R$ cause no exception.
Equations~\eqref{eq:first-zero-count} and \eqref{eq:second-zero-count} contradict $d\geq2$, proving the lemma.
\end{proof}

\begin{figure}[tbp]
\centering
\begin{tikzpicture}[x=1cm,y=1cm,>=Stealth,font=\small]

\begin{scope}[shift={(0,0)},scale=1.08]
  \draw[thin,densely dotted] (0,0) circle (1.45);
  \node at (-0.15,1.73) {$|\lambda|=\rho_*$};

  \coordinate (a) at (50:1.45);
  \coordinate (b) at (-50:1.45);
  \coordinate (c) at (165:1.45);
  \coordinate (d) at (245:1.45);

  \draw[thick] (a)--(c)--(d)--(b)--cycle;

  \begin{scope}[scale=.50]
    \draw[semithick,dashed,dash pattern=on 3pt off 2pt]
      (50:1.45)--(165:1.45)--(245:1.45)--(-50:1.45)--cycle;
  \end{scope}

  \draw[very thick] (a)--(b);

  \fill (a) circle(1.3pt) node[above right=1pt] {$\lambda_2$};
  \fill (b) circle(1.3pt) node[below right=1pt] {$\lambda_1$};
  \fill (c) circle(1.1pt);
  \fill (d) circle(1.1pt);
  \fill (0,0) circle(.9pt);
  \node[below left=1pt] at (0,0) {$0$};

  \draw[->] (.18,0)--(1.95,0) node[above] {$\nu=1$};

  \node at (-0.58,0.38) {$\vartheta K$};
  \node at (-1.08,-0.90) {$K$};

  \node[align=center] at (0,-2.05) {\textit{Frequency plane}};
\end{scope}

\begin{scope}[shift={(4.15,0)}]
  \draw[densely dashed] (0,.40)--(3.20,.40);

  \draw[->,thin] (-.05,-0.92)--(3.18,-0.92) node[right] {$\nu=1$};

  \foreach \x/\lab in {0.45/{k},1.55/{k+1},2.65/{k+2}} {
    \draw[thick] (\x,.40) circle (.28);
    \fill (\x,.40) circle (1pt);
    \node at (\x,-0.13) {$z_{\lab}$};
  }

  \node at (0.45,1.02) {$D_k$};
  \node[align=center] at (1.55,1.45) {$z_k=t_k\nu+z_*$};

  \node[draw,inner sep=4pt,align=center] at (1.55,-1.58)
    {$\displaystyle 1=N_{D_k}(F)=d\,N_{D_k}(h), \qquad d\ge 2$};

  \node[align=center] at (1.55,-2.18) {\textit{Source plane}};
\end{scope}

\end{tikzpicture}

\caption{The exposed-edge obstruction, drawn after inverse rotations
of the frequency and source planes, preserving $\lambda z$,
so that the exposed-edge direction is normalized to $\nu=1$.
The frequencies lie on the circle $|\lambda|=\rho_*$, while the
remainder is confined to the contracted diagram $\vartheta K$.
The zeros $z_k$ of the two-mode sum lie on an affine line in the source
plane. On the fixed-radius discs $D_k$, Rouch\'e's theorem yields
the incompatible count
$1=N_{D_k}(F)=dN_{D_k}(h)$ with $d\ge2$.
The two planes and their scales are independent.}
\label{fig:edge}
\end{figure}

\FloatBarrier

\subsection{Removal of the type-zero remainder}

\begin{proof}[Proof of Theorem~\ref{thm:composition}]
Choose $\lambda\in\Lambda$ and put $b=\lambda/d\neq0$. In the variable $w=bz$, write
\[
 \widetilde h(w)=h(w/b),\qquad V(w)=v(w/b).
\]
By Lemma~\ref{lem:edge}, either $\Lambda=\{\lambda\}$ or $\Lambda=\{\lambda,-\lambda\}$.
Lemma~\ref{lem:poly-diagram}, Lemma~\ref{lem:polygon}, and the scaling rule \eqref{eq:scaling-diagram} give
\begin{equation}\label{eq:segment-diagrams}
 \begin{array}{c|c|c}
 &K(\widetilde h)&\text{a set containing }I(V)\\ \hline
 \Lambda=\{\lambda\}&[0,1]&[0,m_Q]\\
 \Lambda=\{\lambda,-\lambda\}&[-1,1]&[-m_Q,m_Q].
 \end{array}
\end{equation}
The intervals in this display lie on the real axis in the $w$-frequency plane. In both cases,
\begin{equation}\label{eq:normalized-spectral}
 R(\widetilde h(w))=c_+e^{dw}+c_-e^{-dw}+V(w),
\end{equation}
where $c_+\neq0$. Here $c_-=0$ in the one-point case, and $c_-\neq0$ in the two-point case.

Centre the polynomial by setting
\[
 c_0=-\frac{r_{d-1}}{dr_d},\qquad
 \widetilde R(T)=R(c_0+T).
\]
Then
\begin{equation}\label{eq:centered-R}
 \widetilde R(T)=r_dT^d+\sum_{j=0}^{d-2}\widetilde r_jT^j.
\end{equation}
We first prove that, for suitable constants $A\neq0$ and $B$,
\begin{equation}\label{eq:zero-remainder}
 u(w)=\widetilde h(w)-c_0-Ae^w-Be^{-w}
\end{equation}
has exponential type zero. In the one-point case we will take $B=0$.

\smallskip
\noindent\emph{Sectorial inversion.}
Fix $0<\delta<\pi/2$ and consider the closed sector
\[
 S_\delta=\{w\neq0:|\arg w|\leq\pi/2-\delta\}.
\]
For $w\in S_\delta$, put $x=\operatorname{Re}w$, so $x\geq(\sin\delta)|w|$. By \eqref{eq:segment-diagrams} and \eqref{eq:uniform-support},
\begin{equation}\label{eq:V-sector}
 |V(w)|\leq C_\varepsilon\exp\bigl(m_Qx+\varepsilon|w|\bigr)
 \qquad(w\in S_\delta)
\end{equation}
for every $\varepsilon>0$. Since $m_Q<d$, take $\varepsilon=(d-m_Q)\sin\delta/2$ in
\eqref{eq:normalized-spectral}. With $\eta_\delta=(d-m_Q)\sin\delta/2>0$,
we obtain
\begin{equation}\label{eq:sector-leading}
 R(\widetilde h(w))=c_+e^{dw}
 \bigl(1+O_\delta(e^{-\eta_\delta|w|})\bigr)
\end{equation}
uniformly as $|w|\to\infty$ in $S_\delta$.

Set $T(w)=\widetilde h(w)-c_0$.
The modulus of \eqref{eq:sector-leading} tends uniformly to infinity.
Polynomial properness first gives $|T(w)|\to\infty$ uniformly, and leading-term comparison in \eqref{eq:centered-R} then gives
\begin{equation}\label{eq:T-size}
 |T(w)|\asymp e^x\qquad(w\in S_\delta,\ |w|\to\infty).
\end{equation}
Indeed, for large $|T|$, the quantities $|\widetilde R(T)|$ and $|T|^d$ are comparable, whereas \eqref{eq:sector-leading} makes $|\widetilde R(T(w))|$ comparable to $e^{dx}$.
Dividing \eqref{eq:centered-R} by $e^{dw}$, the lower-degree terms are $O(e^{-2x})$ by \eqref{eq:T-size}.
Thus
\[
 \bigl(T(w)e^{-w}\bigr)^d\longrightarrow\frac{c_+}{r_d}
\]
uniformly in the sector.

The roots of $X^d-c_+/r_d$ are distinct and nonzero.
For sufficiently large $|w|$, the quotient $T(w)e^{-w}$ lies in a union of disjoint small discs about these roots.
The truncated sector $S_\delta\cap\{|w|>R_0\}$ is connected, so continuity confines the quotient to one of these discs.
Letting the radii of the discs decrease gives
\begin{equation}\label{eq:chosen-root}
 T(w)e^{-w}\longrightarrow A,\qquad A^d=c_+/r_d,
\end{equation}
uniformly in $S_\delta$, for one fixed root $A\neq0$.
The sectors are nested and all contain the positive real ray, so the root is independent of $\delta$.
This argument selects the root through connectedness; no global branch of an inverse function is required.

Factor
\[
 T^d-A^de^{dw}=(T-Ae^w)
 \sum_{j=0}^{d-1}T^{d-1-j}A^je^{jw}.
\]
By \eqref{eq:chosen-root}, the sum on the right is asymptotic to $dA^{d-1}e^{(d-1)w}$ uniformly in $S_\delta$.
From \eqref{eq:normalized-spectral} and \eqref{eq:centered-R},
\[
 r_d(T^d-A^de^{dw})
 =c_-e^{-dw}+V(w)-\sum_{j=0}^{d-2}\widetilde r_jT^j.
\]
Combining this identity with \eqref{eq:T-size} and \eqref{eq:V-sector} yields
\begin{equation}\label{eq:sector-error}
 |T-Ae^w|\leq C_{\varepsilon,\delta}
 \left(\exp\bigl((m_Q-d+1)x+\varepsilon|w|\bigr)
       +e^{-x}+e^{-(2d-1)x}\right)
\end{equation}
for all sufficiently large $|w|$ in $S_\delta$. The three terms come, respectively, from $V$, the lower powers of $T$, and the opposite exponential mode.

Since $m_Q\leq d-1$ and $\varepsilon$ can be arbitrarily small, \eqref{eq:sector-error} implies that the indicator of $T-Ae^w$ is at most zero on every ray in the open right half-plane.

\smallskip
\noindent\emph{The type-zero conclusion.}
If $c_-=0$, take $B=0$.
On the open left half-plane, $K(\widetilde h)=[0,1]$ already gives indicator at most zero for $\widetilde h$.
The same is true of the constant $c_0$ and of $Ae^w$, so $u$ has indicator at most zero there.

If $c_-\neq0$, repeat the sectorial argument with $w$ replaced by $-w$.
It gives $B\neq0$, with $B^d=c_-/r_d$, such that $\widetilde h-c_0-Be^{-w}$ has indicator at most zero on the open left half-plane.
Subtracting $Ae^w$, which decreases exponentially there, preserves the bound.
On the open right half-plane, subtracting $Be^{-w}$ from \eqref{eq:sector-error} likewise preserves the bound.

In either case, $u$ is of exponential type and has indicator at most zero away from the imaginary axis.
If $u\neq0$, continuity of its indicator gives the same conclusion on the two imaginary directions.
Formula~\eqref{eq:type-diagram} now gives $\tau(u)=0$.
The case $u=0$ is immediate.
Notice that this last step does not require sectorial estimates uniform as $\delta\to0$.

\smallskip
\noindent\emph{Eliminating the remainder.}
In the $w$-variable, the differential equation reads
\begin{equation}\label{eq:scaled-ODE}
 (b^n\DD_w^n-a)R(\widetilde h)=Q(\widetilde h),\qquad a=(db)^n.
\end{equation}
Substitute $\widetilde h=c_0+Ae^w+Be^{-w}+u$.
Each side becomes a finite sum of $e^{jw}$, $j\in\Z$, with entire coefficients of type zero.
Lemma~\ref{lem:zero-independence} therefore allows comparison of the coefficients at each frequency.

The coefficient of $e^{(d-1)w}$ in $R(\widetilde h)=\widetilde R(Ae^w+Be^{-w}+u)$ is
\begin{equation}\label{eq:critical-coefficient}
 dr_dA^{d-1}u.
\end{equation}
To check this, choose $k_+$ factors $Ae^w$, $k_-$ factors $Be^{-w}$, and $k_0$ factors $u$ from the top power. The conditions for frequency $d-1$ are
\[
 k_+-k_-=d-1,\qquad k_++k_-+k_0=d,
\]
so $2k_-+k_0=1$.
Hence $k_-=0$ and $k_0=1$, with multiplicity $d$.
No power of degree at most $d-2$ can reach frequency $d-1$, and the degree-$(d-1)$ term was removed in \eqref{eq:centered-R}.

Let $q_{d-1}$ be the coefficient of $T^{d-1}$ in $Q$, taken as zero if $m_Q<d-1$.
Translation by $c_0$ does not change that coefficient, and the coefficient of $e^{(d-1)w}$ in $Q(\widetilde h)$ is $q_{d-1}A^{d-1}$.
Using
\[
 \DD_w^n\bigl(a_j(w)e^{jw}\bigr)
 =e^{jw}(\DD_w+j)^n a_j(w)
\]
in \eqref{eq:scaled-ODE} and dividing by $dr_dA^{d-1}b^n\neq0$, we obtain
\begin{equation}\label{eq:residual-ODE}
 \left[(\DD_w+d-1)^n-d^n\right]u
 =\frac{q_{d-1}}{dr_db^n}.
\end{equation}
The characteristic roots are
\[
 \rho=d\zeta-(d-1),\qquad \zeta^n=1.
\]
They are simple, and none is zero because $(d-1)/d$ is not a root of unity.
Also $(d-1)^n-d^n\neq0$, so \eqref{eq:residual-ODE} has a constant particular solution $u_0$.
Every solution therefore has the form
\[
 u(w)=u_0+\sum_{\rho}a_\rho e^{\rho w}.
\]
The function $u-u_0$ has type zero.
Apply Lemma~\ref{lem:zero-independence} to this identity with frequency $0$ included: every $a_\rho$ vanishes.
Thus $u$ is constant.
Setting $C=c_0+u_0$ and returning to $z$ proves \eqref{eq:h-classified}.
\end{proof}

\begin{remark}\label{rem:odd-order}
If $n$ is odd, the two-point case cannot occur: $\lambda^n=a$ and $(-\lambda)^n=a$ would imply $a=0$.
In this case the conclusion of Theorem~\ref{thm:composition} has the more specific form $h(z)=Ae^{bz}+C$, with $A,b\neq0$.
\end{remark}

\begin{remark}\label{rem:no-arithmetic-modes}
Theorem~\ref{thm:composition} is purely analytic.
Its proof does not require, or establish, algebraicity of the constants in the spectral decomposition when the Taylor coefficients of the original function are algebraic.
Arithmetic descent will instead be applied to finite polynomial relations in Section~\ref{sec:arithmetic}.
\end{remark}

\section{Meromorphic Briot--Bouquet equations}\label{sec:classification}

We pass from polynomial composition to an arbitrary relation \eqref{eq:BB}.
The following proposition records the two facts needed from the established theory.

\begin{proposition}[Consequences of earlier results]\label{prop:known-input}
Let $f$ be a transcendental meromorphic solution of \eqref{eq:BB} on $\C$.
\begin{enumerate}
\item If $f$ has a pole, then $f\in\W$.
\item If $f$ is entire, then it has order one and finite positive exponential type.
\end{enumerate}
\end{proposition}
\begin{proof}[Derivation]
Theorem~D of \cite{ELN} originates in Eremenko's work
\cite{Eremenko82}; its proof is reproduced in the appendix of
\cite{ELN}. Part~\textup{(i)} is \cite[Theorem~1]{ELN}.
For part~\textup{(ii)}, choose an irreducible factor of the defining polynomial. On its
smooth projective normalization $X$, use the notation of
\cite[Section~3]{ELN}: the coordinate functions $p,q$ satisfy
\[
 p\circ\phi=f^{(n)},\qquad q\circ\phi=f,
\]
where $\phi:\C\to X$ is the holomorphic solution map. Let $I$ be
the pole set of $q$ and, for $\xi\in I$, put
\[
 \kappa(\xi)=\frac{\ord_\xi p}{\ord_\xi q}.
\]
The numerator and denominator are orders of meromorphic functions,
with a pole having negative order. We use the following precise
parts of \cite[Theorem~D]{ELN}, stated here in this notation:
\begin{enumerate}[label=\textup{(\alph*)},start=2]
\item For every $\xi\in I$, either $\kappa(\xi)=1$ or
$\kappa(\xi)=1+n/m$ for an integer $m\geq1$.
\item In the latter case, $\phi^{-1}(\xi)$ is infinite, and each of
its points is a pole of $f$ of order $m$.
\item If $\kappa(\xi)=1$, then $\phi^{-1}(\xi)=\varnothing$.
\end{enumerate}
If $f$ is entire, alternative~\textup{(c)} is impossible: it supplies
actual poles in the source plane, not merely a formal branch at
infinity. Since $q$ is nonconstant on the compact curve $X$, the set
$I$ is nonempty. Hence at least one $\xi\in I$ has
$\kappa(\xi)=1$; in fact all do. Now \cite[Lemma~1]{ELN} gives
order one and normal type.

For an entire function this implies $0<\tau(f)<\infty$ in the
maximum-modulus convention of Section~\ref{sec:borel}. Indeed,
\[
 T(r,f)\leq\log^+M(r,f),\qquad
 \log^+M(r,f)\leq3T(2r,f)+O(1),
\]
and normal type at order one means
$0<\limsup_{r\to\infty}T(r,f)/r<\infty$.
This deduction uses no hypothesis on the multiplicities of the
highest homogeneous factors of $P$.
\end{proof}

\begin{lemma}[Polynomial parametrization]\label{lem:parametrization}
Let $f$ be a transcendental entire solution of \eqref{eq:BB}.
Then either $f\in\W$, or there are a transcendental entire function $h$ of finite positive exponential type and polynomials $R,S\in\C[T]$ of equal positive degree such that
\[
 f=R(h),\qquad f^{(n)}=S(h).
\]
\end{lemma}
\begin{proof}
Choose an irreducible factor of $P$ vanishing on $(f,f^{(n)})$.
Let $X$ be the smooth projective normalization of its plane curve, and let $x,y\in\C(X)$ denote the coordinate functions.
The pair induces a nonconstant holomorphic map
\begin{equation}\label{eq:normalized-map}
 \phi:\C\longrightarrow X,\qquad
 x\circ\phi=f,\quad y\circ\phi=f^{(n)}.
\end{equation}
The meromorphic pair defines a holomorphic map to the projective
plane curve. Since the source is a Riemann surface, this nonconstant
map lifts uniquely by the universal property of normalization.
Neither coordinate on $X$ is constant: a constant $f^{(n)}$ would
make $f$ a polynomial.

If $X$ has genus at least two, lift $\phi$ to the unit disc, its universal cover.
The lift is bounded and entire, hence constant by Liouville's theorem, a contradiction.
If $X$ has genus one, write its universal cover as $\pi:\C\to\C/\Gamma\simeq X$.
Since $\C$ is simply connected, $\phi=\pi\circ H$ for an entire function $H$, necessarily nonconstant.
Picard's theorem says that $H$ omits at most one complex number.
Every fibre of $\pi$ is infinite, so $H$ meets every fibre and $\phi$ is surjective.
But the nonconstant meromorphic function $x$ on the compact curve $X$ has a pole, and its pullback under a surjective nonconstant map has a pole.
This contradicts the fact that $f=x\circ\phi$ is entire.

Therefore $X\simeq\PP^1$.
The union of the pole sets of $x$ and $y$ is nonempty and is omitted by $\phi$, because both functions in \eqref{eq:normalized-map} are entire.
Picard's theorem bounds this union by two points.

Suppose first that it consists of two points.
Choose a coordinate $t$ sending them to $0$ and $\infty$, and put $h=t\circ\phi$.
Then $h$ is entire and zero-free, and $f=R(h)$ for a nonconstant rational function $R$.
The rational-composition identity for the Nevanlinna characteristic \cite[Chapter~1]{Hayman} is
\[
 T(r,R(h))=(\deg R)T(r,h)+O(1),
\]
where $\deg R$ is the degree of the rational map. Proposition~\ref{prop:known-input} gives $T(r,f)=O(r)$, hence $T(r,h)=O(r)$. The standard estimate
\[
 \log M(r,h)\leq3T(2r,h)+O(1)
\]
for an entire function shows that $h$ has finite exponential type.
Hadamard factorization for a zero-free entire function of finite order at most one gives $h(z)=e^{az+b}$; nonconstancy forces $a\neq0$.
Thus $f\in\W$.
This is the two-omitted-point argument also used in \cite[Lemma~2]{ELN}.

Suppose instead that the union has a single point.
Choose a coordinate $t$ on $\PP^1$ whose only pole is there.
Then $x=R(t)$ and $y=S(t)$ are polynomials of positive degree.
Since $\phi$ omits that point, $h=t\circ\phi$ is entire.
It is transcendental, since otherwise $f=R(h)$ would be a polynomial.
Leading-term comparison gives
\[
 \log^+|R(h(z))|=(\deg R)\log^+|h(z)|+O(1).
\]
Taking maxima and using Proposition~\ref{prop:known-input} proves
\[
 0<\tau(h)=\frac{\tau(f)}{\deg R}<\infty.
\]
Finally, Lemmas~\ref{lem:poly-diagram} and \ref{lem:derivative-type} give
\[
 (\deg R)\tau(h)=\tau(f)=\tau(f^{(n)})=(\deg S)\tau(h).
\]
Since $\tau(h)>0$, we obtain $\deg R=\deg S$.
\end{proof}

\begin{remark}[The leading form at infinity]\label{rem:infinity}
In the one-point case of Lemma~\ref{lem:parametrization}, the
coordinate functions generate the function field of the plane curve,
so $\C(R,S)=\C(t)$. Thus $t\mapsto(R(t),S(t))$ is birational
onto that curve, whose projective normalization is $\PP^1$. If $d=\deg R=\deg S$ and $r_d,s_d$ are the
leading coefficients, its projective extension is
\[
 [T:U]\longmapsto
 [\,U^dR(T/U):U^dS(T/U):U^d\,].
\]
It has no base point, is birational, and pulls back a line to a
divisor of degree $d$. Thus the irreducible implicit equation has
total degree $d$. The line at infinity has the single preimage
$[1:0]$, whose image is $[r_d:s_d:0]$. Consequently its highest
homogeneous part is
\[
 c(s_dX-r_dY)^d\qquad(c\neq0).
\]
After scaling, this is exactly \eqref{eq:degenerate-case}.
The two-point case is already a Laurent-exponential case.
The number of points on the normalization should not, however, be
identified with the number of distinct leading linear factors:
different branches may have the same image at infinity.
\end{remark}

\begin{proof}[Proof of Theorem~\ref{thm:classification}]
A rational solution already lies in $\W$, and a transcendental solution with a pole lies in $\W$ by Proposition~\ref{prop:known-input}.
For a transcendental entire solution, apply Lemma~\ref{lem:parametrization}.
Unless membership in $\W$ is already known, we have
\[
 f=R(h),\qquad f^{(n)}=S(h),\qquad \deg R=\deg S=d.
\]
For $d=1$, elimination of $h$ gives $f^{(n)}=af+b$.
For $d\geq2$, Theorem~\ref{thm:composition} makes $f$ a Laurent polynomial in an exponential, hence a member of $\W$.

For the last assertion, suppose that an irreducible $P$ of total degree at least two has a solution outside $\W$.
This solution is nonconstant and satisfies $f^{(n)}=af+b$, so $P(f,af+b)=0$.
A nonzero polynomial in one variable cannot vanish on the infinite range of $f$.
Hence $P(X,aX+b)$ vanishes identically.
Equivalently, $Y-aX-b$ divides $P$, contradicting irreducibility and $\deg P\geq2$.
\end{proof}

\section{Descent and first-order arithmetic rigidity}\label{sec:arithmetic}

We now separate the arithmetic argument from the preceding classification.
The first observation is linear algebra over a Laurent-series field; it does not require the coefficient field to be a number field.

\begin{lemma}[Laurent-series descent]\label{lem:descent}
Let $k\subseteq\C$ be a subfield, and let $F_1,\ldots,F_N\in k((t))$.
If these series are linearly dependent over $\C$, then they are linearly dependent over $k$.
\end{lemma}
\begin{proof}
Choose a common lower bound and write
\[
 F_j=\sum_{m\geq m_0}a_{m,j}t^m,\qquad a_{m,j}\in k.
\]
The rows $(a_{m,1},\ldots,a_{m,N})$ lie in $k^N$.
Select finitely many rows forming a basis, and put them in a matrix $B$.
The rank of $B$ is unchanged by the scalar extension $k\subseteq\C$, since rank is detected by nonzero minors.
A nontrivial complex dependence of the $F_j$ is in the complex kernel of $B$, so this rank is less than $N$.
The kernel over $k$ is therefore also nonzero.
Any nonzero vector in that kernel annihilates every original coefficient row, giving a dependence over $k$.
\end{proof}

\begin{corollary}[Differential-algebraic descent]\label{cor:descent}
Let $k$ be a subfield of $\C$, let $f$ be meromorphic on $\C$, and suppose that $\alpha\in k$ and $f(\alpha+t)\in k((t))$.
Every finite set of distinct differential monomials in $z,f,f',\ldots$ that is linearly dependent over $\C$ is linearly dependent over $k$.
In particular, a nonzero polynomial relation
\[
 P\bigl(z,f,f',\ldots,f^{(m)}\bigr)=0
\]
over $\C$ implies a nonzero such relation over $k$, with support contained in the support of $P$.
\end{corollary}
\begin{proof}
Differentiation and multiplication preserve $k((t))$, and $z=\alpha+t$ has its germ in this field.
Apply Lemma~\ref{lem:descent} to the finitely many monomials in the relation.
A nonzero coefficient vector on distinct formal monomials gives a nonzero polynomial.
The resulting identity holds near $\alpha$ and extends to $\C$ by the identity theorem for meromorphic functions.
\end{proof}

For the rest of the paper, put $k=\Qbar$.
The first-order argument below uses the classical transcendence
theorems of Hermite--Lindemann and Schneider
\cite{Schneider,Lang,Waldschmidt}. Related arithmetic normal forms
for finite-order differential systems of transcendence degree one
are given by Yanchenko \cite[Theorem~1]{Yanchenko25}, under
the assumptions that the solution coordinates are holomorphic and
algebraic at the initial point and that the leading derivative
coefficients and separants are nonzero there.
We give a direct proof on the normalized curve, allowing projective
values and without a separant condition. The ramification identity
used in the proof is the standard local pullback formula for a
differential.
Once the first-order equation is defined over $k$, two algebraic values suffice; no hypothesis on derivatives at these points is needed.

Yanchenko \cite{YanchenkoPairs25} also studies quantitative bounds
for simultaneous algebraic values of two derivative levels of
finite-order meromorphic functions, with the lower derivative outside
$\W$ and with degree, height and magnitude restrictions on the values.
The result below is a qualitative rigidity statement
at two algebraic source points for a first-order equation over $k$.
Further Laurent-exponential normal forms for periodic $E$-functions
over a number field, under the arithmetic coefficient bounds defined
there, are established in \cite[Theorem~1]{Yanchenko26}.

\begin{theorem}[Two-point rigidity]\label{thm:two-point}
Let $f$ be meromorphic on $\C$ and suppose that $P(f,f')=0$ for a nonzero polynomial $P\in k[X,Y]$. If $\alpha,\beta\in k$ are distinct and
\[
 f(\alpha),f(\beta)\in\PP^1(k),
\]
then $f\in k(z)$.
\end{theorem}
\begin{proof}
The constant case follows from either algebraic value.
Suppose that $f$ is nonconstant.
Choose an irreducible factor of $P$ over $k$ vanishing at $(f,f')$.
Since $k$ is algebraically closed, the associated smooth projective normalization $X$ and its coordinate functions $x,y\in k(X)$ are defined over $k$.
As in Lemma~\ref{lem:parametrization}, projective compactification and normalization give a nonconstant holomorphic map
\[
 \phi:\C\longrightarrow X(\C),\qquad
 x\circ\phi=f,\quad y\circ\phi=f'.
\]

\smallskip
\noindent\emph{Algebraic points and the pulled-back differential.}
The nonconstant rational function $x$ defines a finite morphism $X\to\PP^1_k$. Every point of a fibre over $\PP^1(k)$ is defined over $k$: its
residue field is a finite extension of the algebraically closed field
$k$. Hence
\begin{equation}\label{eq:algebraic-curve-points}
 \phi(\alpha),\phi(\beta)\in X(k).
\end{equation}
This includes either of the values $f(\alpha)$ and $f(\beta)$ being infinite.

The function $y$ is not identically zero, since $f$ is nonconstant. Define
\[
 \omega=\frac{dx}{y}.
\]
This is a nonzero rational differential over $k$. On a nonempty open set,
\[
 \phi^*\omega=\frac{df}{f'}=dz,
\]
and the identity extends as an identity of meromorphic differentials everywhere.
To record its local consequence, let $e\geq1$ be the local degree of $\phi$ at $z_0$ and let $m=\ord_{\phi(z_0)}\omega$.
In local parameters, the pullback of a differential of the form $s^m v(s)\,ds$, with $v(0)\neq0$, under $s=a(z-z_0)^e+\cdots$ has order
\[
 em+(e-1)=e(m+1)-1.
\]
Since $dz$ has order zero, we obtain
\begin{equation}\label{eq:local-omega}
 e(m+1)=1.
\end{equation}
Since $e$ and $m+1$ are integers and $e\geq1$, we have $e=1$ and $m=0$.
In particular, $\phi$ is unramified and omits every zero and pole of $\omega$.
This calculation allows zeros of $f'$ and singularities of the plane model; the passage to the normalized curve has removed any need to assume a nonvanishing separant.

Uniformization excludes genus at least two, as in Lemma~\ref{lem:parametrization}. We treat the two remaining genera separately.

\smallskip
\noindent\emph{Genus zero.}
Choose an isomorphism $X\simeq\PP^1_k$.
The support of $\divisor\omega$ is omitted by $\phi$, so Picard's theorem bounds it by two points.
It is nonempty because a rational differential on $\PP^1$ has divisor of degree $-2$.
Its points are defined over $k$.

If the support is a single point, send it to infinity by a coordinate $t$ defined over $k$.
The differential then has the form $\omega=c\,dt$ with $c\in k^\times$: the quotient $\omega/dt$ is a rational function with no zeros or poles.
Since $\phi$ omits infinity, the function $h=t\circ\phi$ is entire.
Pullback gives $ch'=1$.
By \eqref{eq:algebraic-curve-points}, $h(\alpha)\in k$, and therefore
\[
 h(z)=h(\alpha)+\frac{z-\alpha}{c}\in k[z].
\]
As $x\in k(t)$, it follows that $f\in k(z)$.

If the support has two points, send them to $0$ and $\infty$ by a coordinate over $k$. A rational differential with no other zeros or poles has the form
\[
 \omega=ct^m\,dt,\qquad c\in k^\times,\quad m\in\Z.
\]
Now $h=t\circ\phi$ is entire and zero-free, and
\[
 ch^mh'=1.
\]
For $m\neq-1$, this gives
\[
 (h^{m+1})'=\frac{m+1}{c}\neq0.
\]
Even when $m+1<0$, the function $h^{m+1}$ is entire and zero-free.
It cannot be a nonconstant affine function, which has a zero.
Hence $m=-1$, and integration gives
\[
 \frac{h(\beta)}{h(\alpha)}
 =\exp\left(\frac{\beta-\alpha}{c}\right).
\]
The quotient on the left belongs to $k^\times$ by \eqref{eq:algebraic-curve-points}; the exponent is a nonzero algebraic number.
Hermite--Lindemann contradicts this equality; see \cite[Theorem~4]{Waldschmidt}.

\smallskip
\noindent\emph{Genus one.}
Take the algebraic point $O=\phi(\alpha)$ as origin, so that $X$ becomes an elliptic curve over $k$. Choose a Weierstrass model
\[
 E:\quad V^2=4U^3-g_2U-g_3,
 \qquad g_2,g_3\in k,\quad g_2^3-27g_3^2\neq0,
\]
and put $\eta=dU/V$. The complex uniformization $\pi:\C\to E(\C)$ can be normalized by
\[
 \pi(0)=O,\qquad \pi^*\eta=du,
 \qquad \pi(u)=(\wp(u),\wp'(u))
\]
away from the period lattice, where $\wp$ has invariants $g_2,g_3$.

Lift $\phi$ to an entire function $H$ with $\phi=\pi\circ H$ and $H(\alpha)=0$.
The lift is nonconstant.
Each fibre of $\pi$ is infinite, whereas Picard's theorem permits $H$ to omit at most one complex number.
Thus $\phi$ is surjective.
Since $\phi$ omits every zero and pole of $\omega$, there are none: $\omega$ is holomorphic and nowhere zero.

The space of holomorphic differentials on $E$ over $k$ is one-dimensional, so
\[
 \omega=c\eta\qquad(c\in k^\times).
\]
The algebraicity of $c$ also follows directly: $\omega/\eta$ is a constant rational function defined over $k$.
Pullback gives $cH'=1$, and the normalization at $\alpha$ yields
\begin{equation}\label{eq:elliptic-argument}
 H(z)=\frac{z-\alpha}{c},\qquad
 \phi(\beta)=\pi\left(\frac{\beta-\alpha}{c}\right).
\end{equation}
The argument $(\beta-\alpha)/c$ is nonzero and algebraic.
Schneider's theorem states that for algebraic invariants and a nonzero algebraic $u$, the number $u$ is not a pole of $\wp$ and $\wp(u)$ is transcendental; see \cite{Schneider} and \cite[Theorem~20\textup{(1)}]{Waldschmidt}.
Consequently the point on the right of \eqref{eq:elliptic-argument} is not in $E(k)$, contrary to \eqref{eq:algebraic-curve-points}.
The theorem also rules out a nonzero period as the argument, so the projective value $\infty$ introduces no exception.

Only the rational genus-zero case remains, which proves the result.
\end{proof}

\section{Proof of the arithmetic theorems}\label{sec:main-proof}

The affine alternative in Theorem~\ref{thm:classification} is handled directly by the algebraic initial data and Lindemann--Weierstrass.

\Needspace{10\baselineskip}
\begin{lemma}[The affine alternative]\label{lem:affine-arithmetic}
Let $\alpha,\beta\in k$ be distinct, and let $f$ be meromorphic on $\C$ with
\[
 f(\alpha+t)\in k((t)),\qquad f(\beta)\in\PP^1(k).
\]
If $f^{(n)}=af+b$ for some $n\geq1$ and $a,b\in\C$, then $f\in k[z]$.
\end{lemma}
\begin{proof}
If $f$ had a pole of order $m$, then $f^{(n)}$ would have a pole of order $m+n$, whereas $af+b$ would have order at most $m$.
Hence $f$ is entire.
The constant case follows from the germ, so suppose that $f$ is nonconstant.

Apply Corollary~\ref{cor:descent} to the monomials $f^{(n)},f,1$.
The coefficient of $f^{(n)}$ in the resulting nontrivial relation over $k$ cannot vanish, because $1$ and $f$ are linearly independent.
Thus the affine equation has coefficients in $k$.
These coefficients are unique for a nonconstant $f$, so the given $a,b$ themselves belong to $k$.

If $a=0$, repeated integration with the algebraic initial derivatives at $\alpha$ gives a polynomial in $k[z]$.
If $a\neq0$, put $g=f+b/a$.
Then $g^{(n)}=ag$.
Let $\lambda_1,\ldots,\lambda_n$ be the roots of $X^n-a$.
They are distinct, nonzero, and algebraic, and
\begin{equation}\label{eq:affine-expansion}
 g(z)=\sum_{j=1}^n c_je^{\lambda_j(z-\alpha)}.
\end{equation}
For $0\leq\ell<n$, evaluation of the derivatives at $\alpha$ gives
\[
 g^{(\ell)}(\alpha)=\sum_{j=1}^n c_j\lambda_j^\ell\in k.
\]
The coefficient matrix is a Vandermonde matrix with determinant
$\prod_{i<j}(\lambda_j-\lambda_i)\neq0$ in $k$. Hence $c_j\in k$ for every $j$.

Because $f$ is entire, the hypothesis on $f(\beta)$ gives a finite algebraic value, and $g(\beta)\in k$. The numbers
\[
 0,\ \lambda_1(\beta-\alpha),\ldots,\lambda_n(\beta-\alpha)
\]
are distinct and algebraic.
Lindemann--Weierstrass makes their exponentials linearly independent over $k$; see \cite[Theorem~6]{Waldschmidt}.
Evaluating \eqref{eq:affine-expansion} at $\beta$ therefore forces $c_1=\cdots=c_n=0$.
In this case $f=-b/a\in k$, completing the proof.
\end{proof}

\begin{proof}[Proof of Theorem~\ref{thm:germ}]
Apply Theorem~\ref{thm:classification}.
In the affine case, Lemma~\ref{lem:affine-arithmetic} gives the conclusion.
Otherwise $f\in\W$.
Every function in $\W$ satisfies a nonzero first-order relation
\[
 A(f,f')=0,\qquad A\in\C[X,Y].
\]
Indeed, $f$ and $f'$ lie in a common function field of transcendence degree one: $\C(z)$, $\C(e^{az})$, or the field of elliptic functions for a fixed period lattice.

Corollary~\ref{cor:descent} descends such a relation to $k[X,Y]$, using the Laurent germ at $\alpha$.
That germ also gives $f(\alpha)\in\PP^1(k)$.
Theorem~\ref{thm:two-point}, applied at $\alpha$ and $\beta$, now yields $f\in k(z)$.
\end{proof}

\begin{proof}[Proof of Corollary~\ref{cor:allpairs}]
First take $r=0$. If $f$ and $f^{(s)}$ were algebraically dependent,
Theorem~\ref{thm:germ}, with $\alpha=0$ and $\beta=1$, would give
$f\in\Qbar(z)$, contrary to $f\in\MahlerClass$.

Differentiation preserves the algebraicity of every local Laurent
coefficient. It also preserves transcendence for globally meromorphic
functions. Indeed, if $f^{(r)}$ were rational for some $r\geq1$,
then $f$ would have only finitely many poles, since each pole persists
under differentiation. Subtracting their principal parts gives an
entire function $F$ for which $F^{(r)}$ is rational and entire, hence
polynomial. Repeated integration makes $F$ polynomial, so $f$ would
be rational. Therefore $f^{(r)}\in\MahlerClass$. Apply the case
$r=0$ to $f^{(r)}$ and the positive derivative order $s-r$.
\end{proof}

\section{Scope and sharpness}\label{sec:scope}

\subsection{The arithmetic hypotheses}
The second algebraic point in Theorem~\ref{thm:germ} cannot be
omitted: $e^z$ has rational Taylor coefficients at zero and satisfies
$f'=f$, but is not rational. In Theorem~\ref{thm:two-point},
definition of the equation over $\Qbar$ is also necessary. The
function $2^z$ takes algebraic values at zero and one and satisfies
$f'=(\log2)f$, but its germ at zero does not have algebraic
coefficients. The germ hypothesis in Theorem~\ref{thm:germ}
supplies precisely the descent that this example lacks.

Global single-valued meromorphy is essential as well. The branch of
$\sqrt{1+z}$ positive at zero has algebraic Taylor coefficients,
takes an algebraic value at one, and satisfies $2ff'=1$, but is not
rational and does not extend meromorphically to all of $\C$.

\subsection{The field of independence}
Independence over $\C$ is weaker than independence over $\C(z)$.
For example, $u(z)=e^{z^2}$ satisfies $u'=2zu$, while $u,u'$ are
algebraically independent over $\C$. Indeed,
\[
 \C(u,u')=\C(z,u),\qquad z=\frac{u'}{2u},
\]
and $u$ is transcendental over $\C(z)$. This function is not in
$\MahlerClass$; the example only separates the two conclusions.
Likewise, pairwise independence of derivative levels does not give
joint independence of the entire family.

Theorem~\ref{thm:classification} and Theorem~\ref{thm:composition}
use the algebraic curve defined by two functions. Descent itself
has no such restriction, as Corollary~\ref{cor:descent} shows.
Thus the main obstacle to extending the arithmetic argument is a
replacement for the one-dimensional analytic classification.

\section{Further questions and an obstruction}\label{sec:further}

\subsection{Nonautonomous relations and three derivative levels}
A natural next question is whether a function in $\MahlerClass$ can
satisfy a nonzero relation $P(z,f,f^{(n)})=0$, or an autonomous
relation $P(f,f',f'')=0$. The first changes the coefficients along
the source plane; the second need not define a curve. Neither is
covered by Corollary~\ref{cor:allpairs}. Even for a restricted class
of defining polynomials, an analogue of the classification in
Theorem~\ref{thm:classification} would provide a new route to
Mahler's question.

\subsection{Constant-coefficient operators}
It is tempting to replace $\DD^n$ in Theorem~\ref{thm:composition}
by an arbitrary constant-coefficient operator. The following example
shows that this is false, and also prevents an unrestricted extension
of Theorem~\ref{thm:classification} to equations $P(f,Lf)=0$.

\begin{example}\label{ex:general-operator}
Put
\[
 h(z)=e^z+e^{iz},\qquad R(T)=T^2+T,\qquad S(T)=T^2+2T,
\]
and let $L=p(\DD)$, where
\begin{equation}\label{eq:counter-polynomial}
 p(\zeta)=\frac35\zeta^4-(2+2i)\zeta^3
             +6i\zeta^2+(4-4i)\zeta-\frac35.
\end{equation}
One checks that
\[
 p(1)=p(i)=2,\qquad p(2)=p(1+i)=p(2i)=1.
\]
Since the frequencies of $h$ are $1,i$ and those of $h^2$ are
$2,1+i,2i$, it follows that
\[
 L R(h)=S(h).
\]
Yet $h$ is not of the form $Ae^{bz}+Be^{-bz}+C$, by linear
independence of exponentials with distinct frequencies.
More strongly, put $f=R(h)$ and $g=Lf=S(h)$. Then $g-f=h$ and
\begin{equation}\label{eq:general-operator-relation}
 (g-f)^2+g-2f=0.
\end{equation}
The polynomial $(Y-X)^2+Y-2X$ is irreducible: in the coordinates
$T=Y-X,X$ it is $T^2+T-X$. The entire function $f$ is not in
$\W$, because its nonzero frequencies do not lie on one real line,
whereas those of an entire nonpolynomial member of $\W$ lie in
$b\Z$ for some $b\neq0$. Nor does $g=af+b$ hold: comparison at
frequencies $1$ and $2$ would give $a=2$ and $a=1$.
\end{example}

The spectrum of $\zeta^n-a$ consists of simple roots on one circle.
Both properties are used in Lemmas~\ref{lem:polygon} and
\ref{lem:edge}. For more general operators, a meaningful extension
must impose additional spectral conditions or allow a broader
conclusion. Identifying such conditions is a separate problem,
not a consequence of the present contraction argument.

\subsection{Linear differential equations and Mahler's question}
One may also ask whether any $f\in\MahlerClass$ can solve a
nontrivial linear differential equation over $\Qbar(z)$.
This requires an arithmetic specialization argument beyond the
autonomous two-term case. In particular, the two-point conclusion
cannot simply be retained for arbitrary linear equations:
$f(z)=z(z-1)e^z$ has a rational Taylor expansion at zero and
$f(1)=0$, but is transcendental and satisfies
\[
 z(z-1)f'-(z^2+z-1)f=0.
\]
Here the universal condition defining $\MahlerClass$, rather than
two-point data alone, would have to supply the additional rigidity.

\section*{Statements and Declarations}

\subsection*{Funding}
The author was supported by the Brazilian National Council for
Scientific and Technological Development (CNPq), Grant
No.~304467/2023-5.

\subsection*{Competing interests}
The author declares that he has no competing interests.

\subsection*{Data availability}
This article is theoretical; the arguments and examples supporting
its results are contained in the article.

\subsection*{AI Declaration}
OpenAI's ChatGPT was used solely for English-language editing and stylistic corrections. The author assumes full responsibility for the manuscript.


\begin{thebibliography}{99}

\bibitem{ADH}
B. Adamczewski, T. Dreyfus, and C. Hardouin,
Hypertranscendence and linear difference equations,
\emph{J. Amer. Math. Soc.} \textbf{34} (2021), no.~2, 475--503.
\url{https://doi.org/10.1090/jams/960}.

\bibitem{ADHW}
B. Adamczewski, T. Dreyfus, C. Hardouin, and M. Wibmer,
Algebraic independence and linear difference equations,
\emph{J. Eur. Math. Soc.} \textbf{26} (2024), no.~5, 1899--1932.
\url{https://doi.org/10.4171/JEMS/1316}.

\bibitem{BankKaufman}
S.~B. Bank and R.~P. Kaufman,
On Briot--Bouquet differential equations and a question of Einar Hille,
\emph{Math. Z.} \textbf{177} (1981), no.~4, 549--559.
\url{https://doi.org/10.1007/BF01219087}.

\bibitem{BeiseMuller}
H.-P. Beise and J. M\"uller,
Growth of (frequently) hypercyclic functions for differential operators,
\emph{Studia Math.} \textbf{207} (2011), no.~2, 97--115.
\url{https://doi.org/10.4064/sm207-2-1}.

\bibitem{Boas}
R.~P. Boas, Jr.,
\emph{Entire Functions},
Academic Press, New York, 1954.

\bibitem{BriotBouquet}
C. Briot and J.-C. Bouquet,
Int\'egration des \'equations diff\'erentielles au moyen des fonctions elliptiques,
\emph{J. \'Ecole Imp\'er. Polytech.} \textbf{21} (1856), 199--254.

\bibitem{Eremenko82}
A.~E. Eremenko,
Meromorphic solutions of equations of Briot--Bouquet type,
\emph{Teor. Funktsii Funktsional. Anal. i Prilozhen.}
\textbf{38} (1982), 48--56;
English translation, \emph{Amer. Math. Soc. Transl. Ser.~2}
\textbf{133} (1986), 15--23.

\bibitem{EremenkoSurvey}
A.~E. Eremenko,
Meromorphic solutions of algebraic differential equations,
\emph{Russian Math. Surveys} \textbf{37} (1982), no.~4, 61--95;
Russian original, \emph{Uspekhi Mat. Nauk}
\textbf{37} (1982), no.~4, 53--82.
\url{https://doi.org/10.1070/RM1982v037n04ABEH003967}.

\bibitem{EremenkoNote}
A.~E. Eremenko,
Meromorphic solutions of Briot--Bouquet type differential equations,
problem note, March 2011, with addenda through December~14, 2024.
\url{https://www.math.purdue.edu/~eremenko/dvi/de.pdf}
(accessed October~13, 2026).

\bibitem{ELN}
A.~E. Eremenko, L.~W. Liao, and T.~W. Ng,
Meromorphic solutions of higher order Briot--Bouquet differential equations,
\emph{Math. Proc. Cambridge Philos. Soc.}
\textbf{146} (2009), no.~1, 197--206.
\url{https://doi.org/10.1017/S030500410800176X}.

\bibitem{Hayman}
W.~K. Hayman,
\emph{Meromorphic Functions},
Oxford Mathematical Monographs,
Clarendon Press, Oxford, 1964.

\bibitem{Hille}
E. Hille,
Higher order Briot--Bouquet differential equations,
\emph{Ark. Mat.} \textbf{16} (1978), 271--286.
\url{https://doi.org/10.1007/BF02386001}.

\bibitem{HMM}
J. Huang, D. Marques, and M. Mereb,
Algebraic values of transcendental functions at algebraic points,
\emph{Bull. Aust. Math. Soc.}
\textbf{82} (2010), no.~2, 322--327.
\url{https://doi.org/10.1017/S0004972710000316}.

\bibitem{Lang}
S. Lang,
\emph{Introduction to Transcendental Numbers},
Addison--Wesley, Reading, MA, 1966.

\bibitem{Levin}
B.~Ya. Levin,
\emph{Distribution of Zeros of Entire Functions},
revised ed., Translations of Mathematical Monographs, vol.~5,
American Mathematical Society, Providence, RI, 1980.

\bibitem{Mahler}
K. Mahler,
Some suggestions for further research,
\emph{Bull. Austral. Math. Soc.}
\textbf{29} (1984), no.~1, 101--108.
\url{https://doi.org/10.1017/S0004972700021316}.

\bibitem{MerollaCorrection}
Correction (to M.~Merolla, ``A Milloux Type Inequality for Good Linear Operators''),
\emph{Complex Variables Elliptic Equ.},
published online August~11, 2026.
\url{https://doi.org/10.1080/17476933.2026.2713839}.

\bibitem{Picard}
\'E. Picard,
Sur une propri\'et\'e des fonctions uniformes d'une variable
et sur une classe d'\'equations diff\'erentielles,
\emph{C. R. Acad. Sci. Paris}
\textbf{91} (1880), 1058--1061.

\bibitem{Schneider}
T. Schneider,
Arithmetische Untersuchungen elliptischer Integrale,
\emph{Math. Ann.} \textbf{113} (1937), 1--13.
\url{https://doi.org/10.1007/BF01571618}.

\bibitem{Waldschmidt}
M. Waldschmidt,
Elliptic functions and transcendence,
in \emph{Surveys in Number Theory},
Developments in Mathematics, vol.~17,
Springer, New York, 2008, 143--188.
\url{https://doi.org/10.1007/978-0-387-78510-3_7}.

\bibitem{Yanchenko}
A.~Ya. Yanchenko,
One advance in the proof of the conjecture on meromorphic solutions
of Briot--Bouquet type equations,
\emph{Izv. Math.} \textbf{86} (2022), no.~5, 1020--1030.
\url{https://doi.org/10.4213/im9265e}.

\bibitem{Yanchenko24}
A.~Ya. Yanchenko,
Entire solutions of one class of algebraic-differential equations
generalizing the Briot--Bouquet type equations,
\emph{Differ. Equ.} \textbf{60} (2024), no.~8, 1117--1122.
Russian original: \emph{Differ. Uravn.} \textbf{60} (2024), no.~8,
1131--1136; main theorem on p.~1132.
\url{https://doi.org/10.1134/S0012266124080123}.

\bibitem{Yanchenko25}
A.~Ya. Yanchenko,
On a refinement of the Schneider--Lang theorem. II. The arithmetic
of the degenerate case,
\emph{Math. Notes} \textbf{117} (2025), nos.~1--2, 158--167.
\url{https://doi.org/10.1134/S0001434625010158}.

\bibitem{YanchenkoPairs25}
A.~Ya. Yanchenko,
On common algebraic points of a pair of different derivatives of
finite-order meromorphic functions,
\emph{Math. Notes} \textbf{118} (2025), no.~6, 1329--1333.
\url{https://doi.org/10.1134/S0001434625605842}.

\bibitem{Yanchenko26}
A.~Ya. Yanchenko,
On periodic $E$-functions,
\emph{Izv. Math.} \textbf{90} (2026), no.~4, 841--849.
\url{https://doi.org/10.4213/im9744e}.

\end{thebibliography}
\end{document}